\documentclass[12pt]{amsart}
\usepackage{amssymb}

\newtheorem{theorem}{Theorem}[section]
\newtheorem{definition}[theorem]{Definition}
\newtheorem{proposition}[theorem]{Proposition}
\newtheorem{lemma}[theorem]{Lemma}

\begin{document}

\title[Fixed point subfactors]{The planar algebra of a fixed point subfactor}

\author{Teodor Banica}
\address{T.B.: Department of Mathematics, Cergy-Pontoise University, 95000 Cergy-Pontoise, France. {\tt teo.banica@gmail.com}}


\begin{abstract}
We consider inclusions of type $(P\otimes A)^G\subset (P\otimes B)^G$, where $G$ is a compact quantum group of Kac type acting on a ${\rm II}_1$ factor $P$, and on a Markov inclusion of finite dimensional $C^*$-algebras $A\subset B$. In the case $[A,B]=0$, which basically covers all known examples, we show that the planar algebra of such a subfactor is of the form $P(A\subset B)^G$, with $G$ acting in some natural sense on the bipartite graph algebra $P(A\subset B)$.
\end{abstract}

\maketitle

\section*{Introduction}

Many known examples of subfactors appear from quantum groups. This is not surprising, in view of the relation of Jones' work \cite{jo1} with statistical mechanics and quantum field theory \cite{eka}, \cite{tli}, \cite{wit}. Among the quantum group constructions, of particular importance are those of Wenzl \cite{wen}, based on some previous work of Kirillov Jr. \cite{kir}, and of Xu \cite{fxu}, using Drinfeld-Jimbo quantum groups at roots of unity \cite{dri}, \cite{jim}. As remarkable exceptions, we have the Haagerup and Asaeda-Haagerup subfactors \cite{aha}, \cite{haa}.

In this paper we study the class of ``fixed point subfactors'',  introduced in \cite{ba2}. Let $G$ be a compact quantum group in the sense of Woronowicz \cite{wo1}, \cite{wo2}, of Kac type, acting on a ${\rm II}_1$ factor $P$, and on a Markov inclusion of finite dimensional $C^*$-algebras $A\subset B$. Under suitable ergodicity assumptions on the actions, we obtain an inclusion of ${\rm II_1}$ factors $(P\otimes A)^G\subset (P\otimes B)^G$. According to \cite{ba1}, \cite{ba2}, \cite{ba3}, \cite{ba4}, \cite{ba5}, \cite{bbi}, the examples include:
\begin{enumerate}
\item The group-subgroup subfactors.

\item The discrete group subfactors.

\item The projective representation subfactors.

\item The finite index depth 2 subfactors.

\item The subfactors associated to vertex models.

\item The subfactors associated to spin models.

\item Fuss-Catalan subfactors of integer index.

\item Most examples of index 4 subfactors.
\end{enumerate}

The first purpose of this paper is to carefully review the construction of the fixed point subfactors, from \cite{ba2}. Our key observation here will be the fact that, in order for everything to work out properly, one has to make the assumption $[A,B]=0$.

This assumption is of course a bit restrictive, theoretically speaking. However, at the level of concrete examples of subfactors, no one is missed: in the above list, the constructions 1-6 and 8 use $A=\mathbb C$, and the construction 7 uses $B=\mathbb C^n$.

We will prove under this assumption the following general result:

\medskip

\noindent {\bf Theorem.} {\em The planar algebra of $(P\otimes A)^G\subset (P\otimes B)^G$ is $P(A\subset B)^G$, the algebra of $G$-invariant elements of the ``bipartite graph'' planar algebra $P(A\subset B)$.}

\medskip

The proof uses the computation in \cite{ba2} of the higher relative commutants of the subfactor, based on a version of Wassermann's ``invariance principle'' in \cite{was}. With this computation in hand, the problem is to find the correct interpretation of the corresponding planar algebra, as subalgebra of the planar algebra $P(A\subset B)$, constructed by Jones in \cite{jo2}. In the $A=\mathbb C$ case this problem was solved in \cite{ba4}, as a consequence of some more general results, regarding the coactions of non-necessarily Kac algebras. Now in the case of an arbitrary inclusion $A\subset B$, the situation is a priori much more complicated, due to the subtleties with the partition function of $P(A\subset B)$, constructed by Jones in \cite{jo2}. However, the  assumption $[A,B]=0$ simplifies everything, and we get the above result.

The paper is organized as follows: 1 is a preliminary section, in 2-3 we discuss some technical issues, and in 4-5 we state and prove our main results.

\section{Fixed point subfactors}

The fixed point subfactors were introduced in \cite{ba2}, as a unification of several basic constructions of subfactors. These constructions use both groups and group duals, and the natural framework for their unification is that of the compact quantum groups. 

In this paper we will be mainly using unitary compact quantum groups, of Kac type. The axioms here, due to Woronowicz \cite{wo1}, \cite{wo2}, are as follows:

\begin{definition}
A unitary compact quantum group of Kac type is described by a $C^*$-algebra $Z=C(G)$ and a unitary $u\in M_n(Z)$, whose entries generate $Z$, such that:
\begin{enumerate}
\item There exists a morphism $\Delta:Z\to Z\otimes Z$ such that $\Delta(u_{ij})=\Sigma_ku_{ik}\otimes u_{kj}$.

\item There exists a morphism $\varepsilon:Z\to\mathbb C$ such that $\varepsilon(u_{ij})=\delta_{ij}$.

\item There exists a morphism $S:Z\to Z^{opp}$ such that $S(u_{ij})=u_{ji}^*$.
\end{enumerate}
\end{definition}

Here we use the somewhat non-standard letter $Z$ to designate the algebra $C(G)$, because the traditional symbol $A$ will be used for our Markov inclusions $A\subset B$, to be introduced later on. However, as we will soon explain, we are not exactly interested here in $Z=C(G)$, but rather in a certain associated von Neumann algebra $L^\infty(G)$.

The basic example is provided by the compact groups of unitary matrices, $G\subset U_n$. Here $u_{ij}$ are the standard matrix coordinates, $u_{ij}(g)=g_{ij}$, and the maps $\Delta,\varepsilon,S$ as above appear by transposing the usual rule of matrix multiplication $(gh)_{ij}=\Sigma_kg_{ik}h_{kj}$, the unit formula $1_n=(\delta_{ij})$, and the unitary inversion formula $u^{-1}=(u_{ji}^*)$.

Another key example is provided by the duals $\widehat{\Gamma}$ of the finitely generated discrete groups, $\Gamma=<g_1,\ldots,g_n>$. Here we can consider the group algebra $Z=C^*(\Gamma)$, together with the unitary matrix $u=diag(g_1,\ldots,g_n)$, and the maps $\Delta,\varepsilon,S$ as above can be defined on $\Gamma\subset Z$ by $\Delta(g)=g\otimes g$, $\varepsilon(g)=1$, $S(g)=g^{-1}$, then extended to $Z$ by linearity.

In relation with this latter example, observe that the reduced group algebra $Z=C^*_{red}(\Gamma)$ has no counit morphism $\varepsilon:Z\to\mathbb C$, unless $\Gamma$ is amenable. In fact, technically speaking, the Hopf $C^*$-algebras $Z=C(G)$ as defined above are ``full'' in the sense of \cite{wo1}.

In what follows we will be interested in the actions of compact quantum groups on finite von Neumann algebras. Let us recall that a von Neumann algebra $P$ is called ``finite'' when it comes with a faithful positive unital trace $tr:P\to\mathbb C$. We recall also from \cite{wo1} that by performing the GNS construction to $C(G)$ with respect to the Haar integration functional, we obtain a certain von Neumann algebra, denoted $L^\infty(G)$. We have:

\begin{definition}
A coaction of $L^\infty(G)$ on a finite von Neumann algebra $P$ is an injective morphism of von Neumann algebras $\pi:P\to L^\infty(G)\otimes P$ satisfying $(id\otimes\pi)\pi=(\Delta\otimes id)\pi$, $(id\otimes tr)\pi=tr(.)1$ and $\overline{{\mathcal P}}^w=P$, where ${\mathcal P}=\pi^{-1}(C^\infty(G)\otimes_{alg}P)$. The coaction is called:
\begin{enumerate}
\item Ergodic, if the algebra $P^G=\{p\in P|\pi(p)=1\otimes p\}$ reduces to $\mathbb C$.

\item Faithful, if the span of $\{(id\otimes\phi)\pi(P)|\phi\in P_*\}$ is dense in $L^\infty(G)$.

\item Minimal, if it is faithful, and $(P^G)'\cap P=\mathbb C$.
\end{enumerate}
\end{definition}

Let us discuss now a key issue, namely that of taking the tensor product of coactions. As explained in \cite{ba2}, it is impossible to give a fully satisfactory definition here, because of the noncommutativity of $C(G)$. However, it is possible to give a good definition for the fixed point algebra of the ``non-existing'' tensor product of coactions.

Let $C(G')$ be the algebra $C(G)$, taken with the matrix $u^t=(u_{ji})$. Observe that the comultiplication of this algebra is $\Delta'=\Sigma\Delta$, where $\Sigma(a\otimes b)=b\otimes a$. We have:

\begin{definition}
Let $\pi:P\to L^\infty(G')\otimes P$ and $\alpha:A\to L^\infty(G)\otimes A$ be two coactions. The fixed point algebra of their tensor product is
$$(P\otimes A)^G=\{x\in P\otimes A|\Phi(x)=x\otimes 1\}$$
where $\Phi:P\otimes A\to L^\infty(G)\otimes P\otimes A$ is given by $\Phi(p\otimes a)=((S\otimes id)\pi(p))_{12}\alpha(a)_{13}$.
\end{definition}

The basic examples come from the actions of compact groups. Here the linear map $\Phi$ in the statement, which is actually always comultiplicative, is multiplicative as well, and comes by transposition from the usual tensor product of actions.

Some other key examples, where $\Phi$ is not necessarily multiplicative, are discussed in \cite{ba2}. Let us mention here that the above notions are fully understood in the group dual case $G=\widehat{\Gamma}$, and also in the case when $C(G)$ is finite dimensional, and $\alpha=\Delta$.

In what follows we will be interested in the case where $A$ is finite dimensional. Here we know from the general theory that $A$ must be a direct sum of matrix algebras, and that $tr$ appears as a linear combination of the corresponding block traces. However, we won't work at this level of generality, and we will use instead the following key definition:

\begin{definition}
Any finite dimensional $C^*$-algebra $A=\oplus A_i$ can be canonically regarded as a finite von Neumann algebra, in the following way:
\begin{enumerate}
\item The trace is $tr(\oplus x_i)=(\Sigma n_i^2tr(x_i))/(\Sigma n_i^2)$, where $n_i=\dim A_i$.

\item The Hilbert space on which $A$ acts is the $l^2$ space of $tr$.
\end{enumerate}
\end{definition}

The canonical trace has of course some alternative descriptions. For instance it is the unique trace making $\mathbb C\subset A$ a Markov inclusion, or it is the unique trace making $mm^*$ proportional to the identity, where $m:A\otimes A\to A$ is the multiplication. See \cite{ba2}.

In what follows all the finite dimensional algebras, usually denoted $A,B,\ldots$ will be endowed with their canonical traces. We have the following result, from \cite{ba2}:

\begin{proposition}
Let $\pi:P\to L^\infty(G')\otimes P$ be a minimal coaction, and let $\alpha:A\to L^\infty(G)\otimes A$ be a coaction on a finite dimensional algebra. The following are equivalent:
\begin{enumerate}
\item The von Neumann algebra $(P\otimes A)^G$ is a factor.

\item The coaction $\alpha$ is centrally ergodic: $Z(A)\cap A^G=\mathbb C$.
\end{enumerate}
\end{proposition}

Summarizing, we know so far how to construct the algebras $(P\otimes A)^G$, and we know as well when they are factors. So, in order to construct the fixed point subfactors, we just have to consider an inclusion of such factors, coming from an inclusion $A\subset B$.

Before doing so, let us go back to the various requirements on $A,B$, and clarify what are the ``admissible'' inclusions $A\subset B$. The result here, from \cite{ba2}, \cite{ba3}, is as follows:

\begin{proposition}
For an inclusion $A\subset B$, the following are equivalent:
\begin{enumerate}
\item $A\subset B$ commutes with the canonical traces.

\item $A\subset B$ is Markov with respect to the canonical traces.
\end{enumerate}
\end{proposition}

In what follows, we will call ``Markov'' the above type of inclusion. These inclusions are of course of a very special type, for instance they are subject to the  ``triviality of the index'' condition $[B:A]=\dim B/\dim A\in\mathbb N$. We will come back to this subject in sections 2-3 below, with a number of  results regarding such inclusions.

Let us collect now all the above results in a single one, as follows:

\begin{theorem}
Let $\alpha:P\to L^\infty(G')\otimes P$ be a minimal coaction on a finite von Neumann algera, let $A\subset B$ be a Markov inclusion of finite dimensional algebras, and let $\beta:B\to L^\infty(G)\otimes B$ be a coaction which leaves $A$ invariant, and which is centrally ergodic on both $A$ and $B$. Then $(P\otimes A)^G\subset (P\otimes B)^G$ is a subfactor, of index $[B:A]\in\mathbb N$.
\end{theorem}

We refer to \cite{ba2} for details here, and to \cite{ba1}, \cite{ba2}, \cite{ba3}, \cite{ba4}, \cite{ba5}, \cite{bbi} for a number of concrete examples of such subfactors, corresponding to the list in the introduction. 

\section{Markov inclusions}

We recall from section 1 that the inclusions $A\subset B$ which can be used for constructing a fixed point subfactor must be ``Markov'', in the sense that they must commute with the canonical traces. In this section we present an algebraic study of such inclusions.

Let us begin with some basic definitions, from \cite{ghj}:

\begin{definition}
Associated to an inclusion $A\subset B$ of finite dimensional algebras are:
\begin{enumerate}
\item The column vector $(a_i)\in\mathbb N^s$ given by $A=\oplus_{i=1}^sM_{a_i}(\mathbb C)$.

\item The column vector $(b_j)\in\mathbb N^t$ given by $B=\oplus_{j=1}^tM_{b_j}(\mathbb C)$.

\item The inclusion matrix $(m_{ij})\in M_{s\times t}(\mathbb N)$, satisfying $m^ta=b$.
\end{enumerate}
\end{definition}

To be more precise, each minimal idempotent in $M_{a_i}(\mathbb C)$ splits as a sum of minimal idempotents of $B$, and $m_{ij}$ is the number of such idempotents from $M_{b_j}(\mathbb C)$. We have:

\begin{proposition}
For an inclusion $A\subset B$, the following are equivalent:
\begin{enumerate}
\item $A\subset B$ commutes with the canonical traces.

\item We have $mb=ra$, where $r=||b||^2/||a||^2$.
\end{enumerate}
\end{proposition}

\begin{proof}
The weight vectors of the canonical traces of $A,B$ are given by $\tau_i=a_i^2/||a||^2$ and $\tau_j=b_j^2/||b||^2$. By plugging these values into the standard compatibility formula $\tau_i/a_i=\sum_jm_{ij}\tau_j/b_j$, we obtain the condition in the statement.
\end{proof}

We will need as well the following basic facts, from \cite{ghj}:

\begin{definition}
Associated to an inclusion $A\subset B$, with matrix $m\in M_{s\times t}(\mathbb N)$, are:
\begin{enumerate}
\item The Bratteli diagram: this is the bipartite graph $\Gamma$ having as vertices the sets $\{1,\ldots,s\}$ and $\{1,\ldots,t\}$, the number of edges between $i,j$ being $m_{ij}$.

\item The basic construction: this is the inclusion $B\subset A_1$ obtained from $A\subset B$ by reflecting the Bratteli diagram. 

\item The Jones tower: this is the tower of algebras $A\subset B\subset A_1\subset B_1\subset\ldots$ obtained by iterating the basic construction.
\end{enumerate}
\end{definition}

We know that for a Markov inclusion $A\subset B$ we have $m^ta=b$ and $mb=ra$, and so $mm^ta=ra$, which gives an eigenvector for the square matrix $mm^t\in M_s(\mathbb N)$. When this latter matrix has positive entries, by Perron-Frobenius we obtain $||mm^t||=r$.

This equality holds in fact without assumptions on $m$, and we have: 

\begin{theorem}
Let $A\subset B$ be a Markov inclusion, with inclusion matrix $m\in M_{s\times t}(\mathbb N)$.
\begin{enumerate}
\item $r=\dim(B)/\dim(A)$ is an integer.

\item $||m||=||m^t||=\sqrt{r}$.

\item $||\ldots mm^tmm^t\ldots||=r^{k/2}$, for any product of lenght $k$.
\end{enumerate}
\end{theorem}

\begin{proof}
Consider the vectors $a,b$, as in Definition 2.1. We know from definitions and from Proposition 2.2 that we have $b=m^ta$, $mb=ra$, and $r=||b||^2/||a||^2$.

(1) If we construct as above the Jones tower for $A\subset B$, we have, for any $k$:
$$\frac{\dim B_k}{\dim A_k}=\frac{\dim A_k}{\dim B_{k-1}}=r$$

On the other hand, we have as well the following well-known formula:
$$\lim_{k\to\infty}(\dim A_k)^{1/2k}=\lim_{k\to\infty}(\dim B_k)^{1/2k}=||mm^t||$$

By combining these two formulae we obtain $||mm^t||=r$. But from $r\in\mathbb Q$ and $(mm^t)^ka=r^ka$ for any $k\in\mathbb N$, we get $r\in\mathbb N$, and we are done.

(2) This follows from the above equality $||mm^t||=r$, and from the standard equalities $||m||^2=||m^t||^2=||mm^t||$, for any real rectangular matrix $r$.

(3) Let $n$ be the length $k$ word in the statement. First, by applying the norm and by using the formula $||m||=||m^t||=\sqrt{r}$, we obtain the inequality $||n||\leq r^{k/2}$.

For the converse inequality, assume first that $k$ is even. Then $n$ has either $a$ or $b$ as eigenvector (depending on whether $n$ begins with a $m$ or with a $m^t$), in both cases with eigenvalue $r^{k/2}$, and this gives the desired inequality $||n||\geq r^{k/2}$.

Assume now that $k$ is odd, and let $\circ\in\{1,t\}$ be such that $n'=m^\circ n$ is alternating. Since $n'$ has even length, we already know that we have $||n'||=r^{(k+1)/2}$. Together with $||n'||\leq ||m^\circ||\cdot ||n||=\sqrt{r}||n||$, this gives the desired inequality $||n||\geq r^{k/2}$.
\end{proof}

\section{The Jones tower}

Assume that a compact quantum group $G$ acts on a Markov inclusion $A\subset B$, as in section 1. Then $G$ acts on the whole Jones tower for $A\subset B$, with the action being uniquely determined by the fact that it fixes the Jones projections. See \cite{ba2}. We have:

\begin{proposition}
The Jones tower for a fixed point subfactor $(P\otimes A)^G\subset(P\otimes B)^G$ is $(P\otimes A)^G\subset (P\otimes B)^G\subset (P\otimes A_1)^G\subset (P\otimes B_1)^G\subset\ldots$
\end{proposition}

\begin{proof}
The idea is to tensor $P$ with the Jones tower for $A\subset B$, then to remark that the Jones projections for this new tower are invariant under $G$. Together with the abstract characterization of the basic construction in \cite{ghj}, this gives the result. See \cite{ba2}.
\end{proof}

The relative commutants for the above inclusions can be computed as follows:

\begin{proposition}
The relative commutants for the Jones tower $N\subset M\subset N_1\subset M_1\subset\ldots$ associated to fixed point subfactor $(P\otimes A)^G\subset(P\otimes B)^G$ are given by:
\begin{enumerate}
\item $N_s'\cap N_t=(A_s'\cap A_t)^G$.

\item $N_s'\cap M_t=(A_s'\cap B_t)^G$.

\item $M_s'\cap N_t=(B_s'\cap A_t)^G$.

\item $M_s'\cap M_t=(B_s'\cap B_t)^G$.
\end{enumerate}
\end{proposition}

\begin{proof}
As explained in \cite{ba2}, this follows from a suitable quantum group adaptation of Wassermann's ``invariance principle'' in \cite{was}, which basically tells us that ``when computing the higher relative commutants, the part involving the ${\rm II}_1$ factor $P$ dissapears''.
\end{proof}

We use now the fact that for a Markov inclusion, the basic construction and the Jones tower have a particularly simple form. Let us first work out the basic construction:

\begin{proposition}
The basic construction for a Markov inclusion $i:A\subset B$ of index $r\in\mathbb N$ is the inclusion $j:B\subset A_1$ obtained as follows:
\begin{enumerate}
\item $A_1=M_r(\mathbb C)\otimes A$, as an algebra.

\item $j:B\subset A_1$ is given by $mb=ra$.

\item $ji:A\subset A_1$ is given by $(mm^t)a=ra$.
\end{enumerate}
\end{proposition}

\begin{proof}
With notations from the previous section, the weight vector of the algebra $A_1$ appearing from the basic construction is $mb=ra$, and this gives the result.
\end{proof}

For the reminder of this section we fix a Markov inclusion $i:A\subset B$. We have:

\begin{proposition}
The Jones tower $A\subset B\subset A_1\subset B_1\subset\ldots$ associated to a Markov inclusion $i:A\subset B$ is given by:
\begin{enumerate}
\item $A_k=M_r(\mathbb C)^{\otimes k}\otimes A$.

\item $B_k=M_r(\mathbb C)^{\otimes k}\otimes B$.

\item $A_k\subset B_k$ is $id_k\otimes i$.

\item $B_k\subset A_{k+1}$ is $id_k\otimes j$.
\end{enumerate}
\end{proposition}

\begin{proof}
This follows from Proposition 3.3, with the remark that if $i:A\subset B$ is Markov, then so is its basic construction $j:B\subset A_1$.
\end{proof}

Regarding now the relative commutants for this tower, we have here:

\begin{proposition}
The relative commutants for the Jones tower $A\subset B\subset A_1\subset B_1\subset\ldots$ associated to a Markov inclusion $A\subset B$ are given by:
\begin{enumerate}
\item $A_s'\cap A_{s+k}=M_r(\mathbb C)^{\otimes k}\otimes (A'\cap A)$.

\item $A_s'\cap B_{s+k}=M_r(\mathbb C)^{\otimes k}\otimes (A'\cap B)$.

\item $B_s'\cap A_{s+k}=M_r(\mathbb C)^{\otimes k}\otimes (B'\cap A)$.

\item $B_s'\cap B_{s+k}=M_r(\mathbb C)^{\otimes k}\otimes (B'\cap B)$.
\end{enumerate}
\end{proposition}

\begin{proof}
The assertions (1,2,4) follow from Proposition 3.4, and from the general properties of the Markov inclusions. As for the third assertion, observe first that we have:
$$B'\cap A_1=(B'\cap B_1)\cap A_1=(M_r(\mathbb C)\otimes Z(B))\cap (M_r(\mathbb C)\otimes A)=M_r(\mathbb C)\otimes (B'\cap A)$$

This proves the assertion at $s=0,k=1$, and the general case follows from it.
\end{proof}

Observe now that relative commutants in Proposition 3.5 are in general not stable under the action of $G$. We will overcome this problem in the following way:

\begin{definition}
We say that a Markov inclusion $A\subset B$ is abelian if $[A,B]=0$.
\end{definition}

In other words, we are asking for the commutation relation $ab=ba$, for any $a\in A,b\in B$. Note that this is the same as asking that $B$ is an $A$-algebra, $A\subset Z(B)$.

Observe that all inclusions with $A=\mathbb C$ or with $B=\mathbb C^n$ are abelian. This is important for the purposes of the present paper, because, as already pointed out in the introduction, all known examples of fixed point subfactors appear from abelian inclusions. We have:

\begin{proposition}
With the notation $\tilde{B}_k=M_r(\mathbb C)^{\otimes k}\otimes Z(B)$, the relative commutants for the Jones tower $A\subset B\subset A_1\subset B_1\subset\ldots$ of an abelian inclusion are given by:
\begin{enumerate}
\item $A_s'\cap A_{s+k}=A_k$.

\item $A_s'\cap B_{s+k}=B_k$.

\item $B_s'\cap A_{s+k}=A_k$.

\item $B_s'\cap B_{s+k}=\tilde{B}_k$.
\end{enumerate}
\end{proposition}

\begin{proof}
This follows by comparing Proposition 3.4 and Proposition 3.5, and by using the fact that for an abelian inclusion we have $Z(A)=A$, $A'\cap B=B$, $B'\cap A=A$.
\end{proof}

We are now in position of stating and proving the main result in this section. This is an improvement of the previous theoretical results in \cite{ba2}, in the abelian case:

\begin{theorem}
The relative commutants for the Jones tower $N\subset M\subset N_1\subset M_1\subset\ldots$ associated to an abelian fixed point subfactor $(P\otimes A)^G\subset(P\otimes B)^G$ are given by:
\begin{enumerate}
\item $N_s'\cap N_{s+k}=A_k^G$.

\item $N_s'\cap M_{s+k}=B_k^G$.

\item $M_s'\cap N_{s+k}=A_k^G$.

\item $M_s'\cap M_{s+k}=\tilde{B}_k^G$.
\end{enumerate}
\end{theorem}

\begin{proof}
This follows from Proposition 3.2 and Proposition 3.7.
\end{proof}

\section{Planar algebras}

In this section and in the next one we reformulate the general results from the previous section, in terms of Jones' bipartite graph planar algebras \cite{jo2}.

A ``$k$-box'' is a rectangle in the plane, with sides parallel to the real and imaginary axes, having $2k$ marked points on its sides: $k$ on the upper side, and $k$ on the lower side. The points are numbered $1,2,\ldots ,2k$, clockwise starting from top left. We have:

\begin{definition}
A $(k_1,\ldots,k_r,k)$-tangle consists of the following:
\begin{enumerate}
\item Boxes: we have an ``input'' $k_i$-box, one for each $i$, and an ``output'' $k$-box. The input boxes are all disjoint, and are contained in the output box.

\item Strings: all the marked points are paired by strings, which lie inside the output box and outside the input boxes, don't cross, and have to match the parity.

\item Circles: there are also a finite number of closed strings, called circles.
\end{enumerate}
\end{definition}

The tangles can be glued in the obvious way, and the corresponding algebraic structure is the planar operad $\mathcal P$. With this notion in hand, a planar algebra is simply an algebra over $\mathcal P$. Or, in more concrete terms, we have the following definition:

\begin{definition}
A planar algebra is a graded vector space $P=(P_k)$, with multilinear maps $T:P_{k_1}\otimes\ldots\otimes P_{k_r}\to P_k$, one for each tangle, compatible with the gluing operation.
\end{definition}

All this is perhaps a bit too abstract, so let us describe right away a very concrete example. This is the planar algebra of a bipartite graph, constructed by Jones in \cite{jo2}.

Let $\Gamma$ be a bipartite graph, with vertex set $\Gamma_a\cup\Gamma_b$. It is useful to think of $\Gamma$ as being the Bratteli diagram of an inclusion $A\subset B$, in the sense of Definition 2.3.

Our first task is to define the graded vector space $P$. Since the elements of $P$ will be subject to a planar calculus, it is convenient to introduce them ``in boxes'', as follows:

\begin{definition}
Associated to $\Gamma$ is the abstract vector space $P_k$ spanned by the $2k$-loops based at points of $\Gamma_a$. The basis elements of $P_k$ will be denoted 
$$x=\begin{pmatrix}
e_1&e_2&\ldots&e_k\\
e_{2k}&e_{2k-1}&\ldots &e_{k+1}
\end{pmatrix}$$
where $e_1,e_2,\ldots,e_{2k}$ is the sequence of edges of the corresponding $2k$-loop.
\end{definition}

Consider now the adjacency matrix of $\Gamma$, which is of type $M=(^0_{m^t}{\ }^m_0)$. We pick an $M$-eigenvalue $\gamma\neq 0$, and then a $\gamma$-eigenvector $\eta:\Gamma_a\cup\Gamma_b\to\mathbb C-\{0\}$. 

With this data in hand, we have the following construction, due to Jones \cite{jo2}:

\begin{definition}
Associated to any tangle is the multilinear map 
$$T(x_1\otimes\ldots\otimes x_r)=\gamma^c\sum_x\delta(x_1,\ldots ,x_r,x)\prod_m\mu(e_m)^{\pm 1}x$$
where the objects on the right are as follows:
\begin{enumerate}
\item The sum is over the basis of $P_k$, and $c$ is the number of circles of $T$.

\item $\delta=1$ if all strings of $T$ join pairs of identical edges, and $\delta=0$ if not.

\item The product is over all local maxima and minima of the strings of $T$. 

\item $e_m$ is the edge of $\Gamma$ labelling the string passing through $m$ (when $\delta=1$).

\item $\mu(e)=\sqrt{\eta(e_f)/\eta(e_i)}$, where $e_i,e_f$ are the initial and final vertex of $e$.

\item The $\pm$ sign is $+$ for a local maximum, and $-$ for a local minimum.
\end{enumerate}
\end{definition}

In other words, we plug the loops $x_1,\ldots,x_r$ into the input boxes of $T$, and then we construct the ``output'': this is the sum of all loops $x$ satisfying the compatibility condition $\delta=1$, altered by certain normalization factors, coming from the eigenvector $\eta$.

Let us work out now the precise formula of the action, for 6 carefully chosen tangles, which are of key importance for the considerations to follow. This study will be useful as well as an introduction to Jones' result in \cite{jo2}, stating that $P$ is a planar algebra:

\begin{definition}
We have the following examples of tangles:
\begin{enumerate}
\item Identity $1_k$: the $(k,k)$-tangle having $2k$ vertical strings.

\item Multiplication $M_k$: the $(k,k,k)$-tangle having $3k$ vertical strings.

\item Inclusion $I_k$: the $(k,k+1)$-tangle like $1_k$, with an extra string at right.

\item Shift $J_k$: the $(k,k+2)$-tangle like $1_k$, with two extra strings at left.

\item Expectation $U_k$: the $(k+1,k)$-tangle like $1_k$, with a curved string at right.

\item Jones projection $E_k$: the $(k+2)$-tangle having two semicircles at right.
\end{enumerate}
\end{definition}

Let us look first at the identity $1_k$. Since the solutions of $\delta(x,y)=1$ are the pairs of the form $(x,x)$, this tangle acts by the identity:
$$1_k\begin{pmatrix}f_1&\ldots&f_k\\ e_1&\ldots&e_k\end{pmatrix}=
\begin{pmatrix}f_1&\ldots&f_k\\ e_1&\ldots&e_k\end{pmatrix}$$

A similar argument applies to the multiplication $M_k$, which acts as follows:
$$M_k\left( 
\begin{pmatrix}f_1&\ldots&f_k\\ e_1&\ldots&e_k\end{pmatrix}
\otimes\begin{pmatrix}h_1&\ldots&h_k\\ g_1&\ldots&g_k\end{pmatrix}
\right)=
\delta_{f_1g_1}\ldots\delta_{f_kg_k}
\begin{pmatrix}h_1&\ldots&h_k\\ e_1&\ldots& e_k\end{pmatrix}$$

Regarding now the inclusion $I_k$, the solutions of $\delta(x_0,x)=1$ being the elements $x$ obtained from $x_0$ by adding to the right a vector of the form $(^g_g)$, we have:
$$I_k\begin{pmatrix}f_1&\ldots&f_k\\ e_1&\ldots&e_k\end{pmatrix}=
\sum_g\begin{pmatrix}f_1&\ldots&f_k&g\\ e_1&\ldots&e_k&g\end{pmatrix}$$

The same method applies to the shift $J_k$, whose action is given by:
$$J_k\begin{pmatrix}f_1&\ldots&f_k\\ e_1&\ldots&e_k\end{pmatrix}=
\sum_{gh}\begin{pmatrix}g&h&f_1&\ldots&f_k\\ g&h&e_1&\ldots&e_k\end{pmatrix}$$

Let us record some partial conclusions, coming from the above formulae:

\begin{proposition}
The graded vector space $P=(P_k)$ becomes a graded algebra, with the multiplication $xy=M_k(x\otimes y)$ on each $P_k$, and with the above inclusion maps $I_k$. The shift $J_k$ acts as an injective morphism of algebras $P_k\to P_{k+2}$.
\end{proposition}

\begin{proof}
The fact that the multiplication is indeed associative follows from its above formula, which is nothing but a generalization of the usual matrix multiplication. The assertions about the inclusions and shifts follow as well by using their above explicit formula.
\end{proof}

Let us go back now to the remaining tangles in Definition 4.5. The usual method applies to the expectation tangle $U_k$, which acts with a spin factor, as follows:
$$U_k\begin{pmatrix}f_1&\ldots&f_k&h\\ e_1&\ldots&e_k&g\end{pmatrix}=
\delta_{gh}\mu(g)^2\begin{pmatrix}f_1&\ldots&f_k\\ e_1&\ldots&e_k\end{pmatrix}$$

As for the Jones projection $E_k$, this tangle has no input box, so we can only apply it to the unit of $\mathbb C$. And when doing so, we obtain the following element:
$$E_k(1)=\sum_{egh}\mu(g)\mu(h)\begin{pmatrix}e_1&\ldots &e_k&h&h\\ e_1&\ldots&e_k&g&g\end{pmatrix}$$

Once again, let us record now some partial conclusions, coming from these formulae:

\begin{proposition}
The elements $e_k=\gamma^{-1}E_k(1)$ are projections, and define a representation of the Temperley-Lieb algebra $TL(\gamma)\to P$. The maps $U_k$ are bimodule morphisms with respect to $I_k$, and their composition is the canonical trace on the image of $TL(\gamma)$.
\end{proposition}

\begin{proof}
The proof of all the assertions is standard, by using the fact that $\eta$ is a $\gamma$-eigenvector of the adjacency matrix. Note that the statement itself is just a generalization of the usual Temperley-Lieb algebra representation on tensors, from \cite{jo4}.
\end{proof}

A careful look at the above computations shows that the following phenomenon appears: the gluing of tangles always corresponds to the composition of multilinear maps.

This phenomenon holds in fact in full generality: the graded linear space $P=(P_k)$, together with the action of the planar tangles given in Definition 4.4, is a planar algebra. We refer to Jones' paper \cite{jo2} for full details regarding this result.

Let us go back now to the Markov inclusions $A\subset B$, as in section 3. We have here the following result from Jones' paper \cite{jo2}, under the assumption that $A$ is abelian:

\begin{theorem}
The planar algebra associated to the graph of $A\subset B$, with eigenvalue $\gamma=\sqrt{r}$ and eigenvector $\eta(i)=a_i/\sqrt{\dim A}$, $\eta(j)=b_j/\sqrt{\dim B}$, is as follows:
\begin{enumerate}
\item The graded algebra structure is given by $P_{2k}=A'\cap A_k$, $P_{2k+1}=A'\cap B_k$.

\item The elements $e_k$ are the Jones projections for $A\subset B\subset A_1\subset B_1\subset\ldots$

\item The expectation and shift are given by the above formulae.
\end{enumerate}
\end{theorem}

\begin{proof}
As a first observation, $\eta$ is indeed a $\gamma$-eigenvector for the adjacency matrix of the graph. Indeed, by using the formulae $m^ta=b$, $mb=ra$, $\sqrt{r}=||b||/||a||$, we get:
$$\begin{pmatrix}0&m\\ m^t&0\end{pmatrix}\begin{pmatrix}a/||a||\\ b/||b||\end{pmatrix}=\begin{pmatrix}\gamma^2a/||b||\\ b/||a||\end{pmatrix}=\gamma\begin{pmatrix}\gamma a/||b||\\ b/\gamma||a||\end{pmatrix}=\gamma\begin{pmatrix}a/||a||\\ b/||b||\end{pmatrix}$$

Since the algebra $A$ was supposed abelian, the Jones tower algebras $A_k,B_k$ are simply the span of the $4k$-paths, respectively $4k+2$-paths on $\Gamma$, starting at points of $\Gamma_a$. With this description in hand, when taking commutants with $A$ we have to just have to restrict attention from paths to loops, and we obtain the above spaces $P_{2k},P_{2k+1}$. See \cite{jo2}.
\end{proof}

\section{Invariant algebras}

In this section we state and prove the main result. Let us begin with a reformulation of Theorem 4.8, in the particular case of the inclusions satisfying $[A,B]=0$:

\begin{proposition}
The ``bipartite graph'' planar algebra $P(A\subset B)$ associated to an abelian inclusion $A\subset B$ can be described as follows:
\begin{enumerate}
\item As a graded algebra, this is the Jones tower $A\subset B\subset A_1\subset B_1\subset\ldots$

\item The Jones projections and expectations are the usual ones for this tower.

\item The shifts correspond to the canonical identifications $A_1'\cap P_{k+2}=P_k$.
\end{enumerate}
\end{proposition}

\begin{proof}
The first assertion is just a reformulation of Theorem 4.8 in the abelian case, by using the identifications $A'\cap A_k=A_k$ and $A'\cap B_k=B_k$ coming from Proposition 3.7.

The assertion on Jones projections follows as well from Theorem 4.8, and the assertion on expectations follows from the fact that their composition is the usual trace.

Regarding now the third assertion, let us recall first from Proposition 3.7 that we have indeed identifications $A_1'\cap A_{k+1}=A_k$ and $A_1'\cap B_{k+1}=B_k$. By using the path model for these algebras, as in the proof of Theorem 4.8, we obtain the result.
\end{proof}
 
In order to formulate now our main result, we will need a few abstract notions, coming from the results in \cite{ba4}. Let $G$ be a compact quantum group, as in section 1. We have:

\begin{definition}
Let $P_1,P_2$ be two finite dimensional algebras, coming with coactions $\alpha_i:P_i\to L^\infty(G)\otimes P_i$, and let $T:P_1\to P_2$ be a linear map.
\begin{enumerate}
\item We say that $T$ is $G$-equivariant if $(id\otimes T)\alpha_1=\alpha_2T$.

\item We say that $T$ is weakly $G$-equivariant if $T(P_1^G)\subset P_2^G$.
\end{enumerate}
\end{definition}

Consider now a planar algebra $P=(P_k)$. The annular category over $P$ is the collection of maps $T:P_k\to P_l$ coming from the ``annular'' tangles, having at most one input box. These maps form sets $Hom(k,l)$, and these sets form a category \cite{jo3}. We have:

\begin{definition}
A coaction of $L^\infty(G)$ on a planar algebra $P$ is a graded algebra coaction $\alpha:P\to L^\infty(G)\otimes P$, such that the annular tangles are weakly $G$-equivariant.
\end{definition}

This definition might seem a bit clumsy, and indeed, it is so: due to a relative lack of concrete examples, this is the best notion of coaction that we have so far. In fact, as it will be shown below, the examples are basically those coming from actions of compact quantum groups on Markov inclusions $A\subset B$, under the assumption $[A,B]=0$.

For the moment, however, let us remain at the generality level of Definition 5.3:

\begin{proposition}
If $G$ acts on on a planar algebra $P$, then $P^G$ is a planar algebra.
\end{proposition}

\begin{proof}
The weak equivariance condition tells us that the annular category is contained in the suboperad $\mathcal P'\subset\mathcal P$ consisting of tangles which leave invariant $P^G$. On the other hand the multiplicativity of $\alpha$ gives $M_k\in\mathcal P'$, for any $k$. Now since $\mathcal P$ is generated by multiplications and annular tangles, we get $\mathcal P'=\mathcal P$, and we are done.
\end{proof}

Let us go back now to the abelian inclusions. We have the following key lemma: 

\begin{lemma}
If $G$ acts on an abelian inclusion $A\subset B$, the canonical extension of this coaction to the Jones tower is a coaction of $G$ on the planar algebra $P(A\subset B)$.
\end{lemma}

\begin{proof}
We know from Proposition 5.1 that, as a graded algebra, $P=P(A\subset B)$ coincides with the Jones tower $A\subset B\subset A_1\subset B_1\subset\ldots$ Thus the canonical Jones tower coaction in the statement can be regarded as a graded coaction $\alpha:P\to L^\infty(G)\otimes P$.  

We have to prove that the annular tangles are weakly equivariant. For this purpose, we use a standard method, from \cite{ba4}. First, since the annular category is generated by $I_k,E_k,U_k,J_k$, we just have to prove that these 4 particular tangles are weakly equivariant. Now since $I_k,E_k,U_k$ are plainly equivariant, by construction of the coaction of $G$ on the Jones tower, it remains to prove that the shift $J_k$ is weakly equivariant.

For this purpose, we use the results in section 3. We know from Theorem 3.8 that the image of the fixed point subfactor shift $J_k'$ is formed by the $G$-invariant elements of the relative commutant $A_1'\cap P_{k+2}=P_k$. Now since this commutant is the image of the planar shift $J_k$ from Proposition 5.1, we have $Im(J_k)=Im(J_k')$, and this gives the result.
\end{proof}

With this lemma in hand, we can now prove:

\begin{proposition}
Assume that $G$ acts on an abelian inclusion $A\subset B$. Then the graded vector space of fixed points $P(A\subset B)^G$ is a planar subalgebra of $P(A\subset B)$.
\end{proposition}

\begin{proof}
This follows from Proposition 5.4 and Lemma 5.5.
\end{proof}

We are now in position of stating and proving our main result:

\begin{theorem}
In the abelian case, the planar algebra of $(P\otimes A)^G\subset (P\otimes B)^G$ is the fixed point algebra $P(A\subset B)^G$ of the bipartite graph algebra $P(A\subset B)$.
\end{theorem}

\begin{proof}
Let $P=P(A\subset B)$, and let $Q$ be the planar algebra of the fixed point subfactor. We know from Theorem 3.8 that we have an equality of graded algebras $Q=P^G$.

It remains to prove that the planar algebra structure on $Q$ coming from the fixed point subfactor agrees with the planar algebra structure of $P$, coming from Proposition 5.1.

Since $\mathcal P$ is generated by the annular category $\mathcal A$ and by the multiplication tangles $M_k$, we just have to check that the annular tangles agree on $P,Q$. Moreover, since $\mathcal A$ is generated by $I_k,E_k,U_k,J_k$, we just have to check that these tangles agree on $P,Q$.

We know that $Q\subset P$ is an inclusion of graded algebras, that all the Jones projections  for $P$ are contained in $Q$, and that the conditional expectations agree. Thus the tangles $I_k,E_k,U_k$ agree on $P,Q$, and the only verification left is that for the shift $J_k$.

Now by using either the axioms of Popa in \cite{pop}, or the construction of Jones in \cite{jo4}, it is enough to show that the image of the subfactor shift $J_k'$ coincides with that of the planar shift $J_k$. But this follows from the argument in the proof of Lemma 5.5.
\end{proof}

There are several questions arising from the present work, the first of which would be the abstract characterization of the planar algebras that we can obtain in this way. This is part of a more general question, regarding the structure of the subfactors having integer index. To our knowledge, nothing much is known here, besides the fact, going back to \cite{ppo}, that the Pimsner-Popa basis appears in a ``clean'' way, without need to complete. However, exploiting this old and well-known fact is a difficult problem.

\end{document}